\documentclass{amsart}

\newtheorem{theorem}{Theorem}[section]
\newtheorem{lemma}[theorem]{Lemma}

\theoremstyle{definition}

\theoremstyle{remark}

\numberwithin{equation}{section}

\begin{document}

\title[Asymptotic behavior for the radial eigenvalues of $p$-Laplacian]{Asymptotic behavior for the radial eigenvalues of $p$-Laplacian in certain annular domains}

\author{Anderson L. A. de Araujo}
\address{Departamento de Matem\'atica, Universidade Federal de Vi\c cosa, 36570-900, Vi\c cosa (MG), Brazil}
\email{anderson.araujo@ufv.br}






\begin{abstract}
In this paper we prove an asymptotic behavior for the radial eigenvalues to the Dirichlet $p$-Laplacian problem $-\Delta_p\,u = \lambda\,|u|^{p-2}u$ in $\Omega$, $u=0$ on $\partial\Omega$, where $\Omega$ is an annular domain $\Omega=\Omega_{R,\overline{R}}$ in $\mathbb{R}^N$. 
\end{abstract}

\maketitle

\section{Introduction}
\label{intro}
This paper investigates an asymptotic behavior for the radial eigenvalues $\lambda_k=\lambda_k(R,\overline{R})$ (when $0<\overline{R}=R+1$ and $R \to +\infty$) to the following eigenvalue problem
\begin{equation}
\label{P5}
\left\{
\begin{array}{lll}
\displaystyle -\Delta_p\,u = \lambda\,|u|^{p-2}u & \textup{in} &\Omega_{R,\overline{R}},\\
\displaystyle u=0 &\textup{on}& \partial\Omega_{R,\overline{R}},
\end{array}
\right.
\end{equation} 
where $\Omega_{R,\overline{R}}= \{x \in \mathbb{R}^N: R<|x|<\overline{R}\}$, with $0<R<\overline{R}$ constants in $\mathbb{R}$, is the annular domain, and we suppose that 
\[1<p\leq\,N.\] 
In particular, when $p=2$, we obtain the Dirichlet Laplacian problem
\begin{equation}
\label{P5.1}
\left\{
\begin{array}{rll}
\displaystyle -\Delta\,u = \lambda\,u & \textup{in} &\Omega_{R,\overline{R}},\\
\displaystyle u=0 &\textup{on}& \partial\Omega_{R,\overline{R}}.
\end{array}
\right.
\end{equation}
Since we are interested only in the radial eigenvalues of \eqref{P5}, we can
rewrite \eqref{P5} as the following $1$-dimensional eigenvalue problem
\begin{equation}
\label{P5.2}
\left\{
\begin{array}{lll}
\displaystyle (r^{N-1}|u'(r)|^{p-2}u'(r))'+\lambda\,r^{N-1}|u(r)|^{p-2}u(r) = 0 & \textup{in} &(R,\bar{R}),\\
\displaystyle u(R)=u(\bar{R})=0. && 
\end{array}
\right.
\end{equation}
We remark that for every $1\leq k \in \mathbb{N}$, if we denote by $\lambda_k$ the $k$-th eigenvalue of
\eqref{P5.2} and by $\lambda_k^{\mbox{rad}}$ the $k$-th radial eigenvalue of \eqref{P5},
\[
\lambda_k=\lambda_k^{\mbox{rad}}.
\]
In order to study the solution of \eqref{P5}, one can make a standard change of variables, see for example \cite{Araujo,LYang}. 

If $N>p$, let $t=-\frac{A}{r^{(N-p)/(p-1))}}+B$ and $v(t)=u(r)$, where 
\[A=\frac{(R\overline{R})^{\frac{N-p}{p-1}}}{\overline{R}^{\frac{N-p}{p-1}} - R^{\frac{N-p}{p-1}}}\,\,\,\mbox{and}\,\,\, B=\frac{\overline{R}^{\frac{N-p}{p-1}}}{\overline{R}^{\frac{N-p}{p-1}} - R^{\frac{N-p}{p-1}}},\] 
 then the problem \eqref{P5.2} (hence \eqref{P5}) transforms into the boundary value problem for the nonautonomous ODE
\begin{equation}
\label{P6}
\left\{
\begin{array}{lll}
\displaystyle (|v'(t)|^{p-2}v'(t))'+\lambda\,q(t)|v(t)|^{p-2}v(t) = 0 & \textup{in} &(0,1),\\
\displaystyle v(0)=v(1)=0. && 
\end{array}
\right.
\end{equation}
where 
\begin{equation}\label{q}
	q(t):=q_{R,\overline{R}}(t)=\left(\frac{p-1}{N-p}\right)^{p}\frac{A^{\frac{(p-1)p}{N-p}}}{(B-t)^{\frac{p(N-1)}{N-p}}}.
\end{equation}

In the case $p=N$, one sets $r=R(\frac{\overline{R}}{R})^t$ and $v(t)=u(r)$, obtaining again the problem \eqref{P6}, now with
\begin{equation}\label{q.2}
	q(t):=q_{R,\overline{R}}(t)=\left[R\left(\frac{\overline{R}}{R}\right)^t\ln\frac{\overline{R}}{R}\right]^{p}.
\end{equation}

Let $\lambda_{k}$ be the $k$-th eigenvalue of \eqref{P6} and let $\varphi_{k}$ be an eigenfunction corresponding to $\lambda_{k}$. Since $q$ satisfies
\[q \in C^{1}([0, 1]), \,\,\, q(t)>0\,\,\, \mbox{for} \,\,\, 0\leq t\leq 1.\]
It is known that
\[0<\lambda_{1}<\lambda_{2}<\cdots<\lambda_{k}<\lambda_{k+1}<\cdots , \,\,\, \lim_{k \to \infty}\lambda_{k}=\infty,\]
and that $\varphi_{k}$ has exactly $k-1$ zeros in $(0, 1)$, see \cite{DR,KN}.

Motivated by the work of Zhang \cite{Zhang}, whose purpose was to compute an estimate for eigenvalues of the Dirichlet $p$-Laplacian \eqref{P6}, $p>1$, we propose to prove an asymptotic behavior for the $\lambda_k(R,\overline{R})$ in the form:
\[\lim_{R \to +\infty} \lambda_k(R,R+1).\]  
The following estimate is known, see Zhang \cite[Remark 2.1]{Zhang}. We suppose that
\begin{equation}\label{pi_p}
	\pi_p=\frac{2\pi(p-1)^{1/p}}{p\sin(\pi/p)},
\end{equation}
\[ 
\bar{q}_-(R,\overline{R}) = \int_0^1\min\{1,q(t)\}dt
\]
and
\[
\bar{q}_+(R,\overline{R}) = \int_0^1\max\{1,q(t)\}dt.
\] 
In \cite{Zhang} (with $T=1$), the author proved the inequality
\begin{equation}\label{estlamb}
	\left(\frac{k\pi_p}{\bar{q}_+(R,\overline{R})}\right)^p\leq \lambda_k(R,\overline{R}) \leq \left(\frac{k\pi_p}{\bar{q}_-(R,\overline{R})}\right)^p.
\end{equation}

In S.S. Lin \cite[Lemma A.1]{Lin}, the author proves an asymptotic behavior for all eigenvalues (that is, radial and non-radial eigenvalues), of Dirichlet problem \eqref{P5.1}, which is the following result.
\begin{lemma}[\cite{Lin}]
Let $\lambda_{k,j}(R,R+1)$ be the $j$-th eigenvalue of 
\[
\phi''+\frac{N-1}{r}\phi'-\frac{\alpha_k}{r^2}\phi=-\lambda_{k,j}(R,R+1)\phi, \mbox{ in } (R,R+1),
\]
\[
\phi(R)=\phi(R+1)=0,
\]
where $\alpha_k=k(k+N-2)$, and let $\lambda_j=j^2\pi^2$ be the $j$-th eigenvalue of
\[
\phi''=-\lambda_{j}\phi, \mbox{ in } (0,1),
\]
\[
\phi(0)=\phi(1)=0,
\]
where $k=0,1,2\cdots$ and $j=1,2,3,\cdots$. Then
\[
\lim\limits_{R\to +\infty}\lambda_{k,j}(R,R+1)=\lambda_j.
\]
\end{lemma}
In this paper, we will prove a generalization of the results of \cite{Lin}, for the radial eigenvalues of $p$-Laplacian, in the cases when $p=r+1$ and $N=2r+1$, with $r \in \mathbb{N}$; $p=N$ and $p=2<3\leq N$. The last case is another proof of Lin's result for the radial eigenvalues. It is noteworthy that it is not yet known a characterization for all eigenvalues of \eqref{P5}. The results of \cite{Lin} are consequence of some results of Bessel functions, while in our paper we use a totally different approach following the work of Zhang \cite{Zhang}.

We state the main result.

\begin{theorem}\label{T2}
$\textbf{(i)}$ Suppose that $N=p$. Then
\begin{equation}\label{fundlim3.2}
	\lim_{R \to +\infty}\bar{q}_+(R,R+1)=\lim_{R \to +\infty}\bar{q}_-(R,R+1)=1.
\end{equation}
In particular, by \eqref{estlamb}
\begin{equation}\label{fundlim4.2}
	\lim_{R \to +\infty} \lambda_k(R,R+1) = \left(k\pi_p\right)^p,
\end{equation}
that is, the eigenvalues $\lambda_k(R,R+1)$ converge asymptotically to the eigenvalues of the problem
  \begin{equation}
\label{P7.2.1}
\left\{
\begin{array}{lll}
\displaystyle (|v'(t)|^{p-2}v'(t))'+\lambda\,|v(t)|^{p-2}v(t) = 0 & \textup{in} &(0,1),\\
\displaystyle v(0)=v(1)=0. && 
\end{array}
\right.
\end{equation}
$\textbf{(ii)}$ Suppose that $p=2$ and $N\geq 3$. Then,
\begin{equation}\label{fundlim}
	\lim_{R \to +\infty}\bar{q}_+(R,R+1)=\lim_{R \to +\infty}\bar{q}_-(R,R+1)=1.
\end{equation}
In particular, by \eqref{estlamb}
\begin{equation}\label{fundlim2}
	\lim_{R \to +\infty} \lambda_k(R,R+1) = \left(k\pi\right)^2,
\end{equation}
that is, the eigenvalues $\lambda_k(R,R+1)$ converge asymptotically to the eigenvalues of the problem
  \begin{equation}
\label{P7}
\left\{
\begin{array}{lll}
\displaystyle v''(t)+\lambda\,v(t) = 0 & \textup{in} &(0,1),\\
\displaystyle v(0)=v(1)=0. && 
\end{array}
\right.
\end{equation}
\end{theorem}

\section{Proofs }\label{section 2}

\begin{proof}[\textbf{Proof of Theorem \ref{T2}}] 
\textbf{(i):}
By \eqref{q.2}, we have
\begin{equation}\label{q2.0}
	q(t)=\left[R\left(\frac{\overline{R}}{R}\right)^t\ln\frac{\overline{R}}{R}\right]^{p}.
\end{equation}
Therefore,
\[
q'(t)=p\left[R\left(\frac{\overline{R}}{R}\right)^t\ln\frac{\overline{R}}{R}\right]^{p-1}R\left(\frac{\overline{R}}{R}\right)^t\left(\ln\frac{\overline{R}}{R}\right)^2>0,
\]
hence, the function $q$ is increasing at $t$ and
\[
0<q(0)\leq q(t) \leq q(1),\,\,\ \forall \,\, t \in [0,1],
\]
that is,
\[
\left[R\ln\frac{\overline{R}}{R}\right]^p\leq q(t) \leq \left[\overline{R}\ln\frac{\overline{R}}{R}\right]^p,\,\,\ \forall \,\, t \in [0,1].
\]

If $R>0$ and $\overline{R}=R+1$, we have
\[
q(0)=\left[R\ln\left(\frac{R+1}{R}\right)\right]^p=\left[\ln\left(1+\frac{1}{R}\right)^R\right]^p
\]
and
\[
q(1)=\left[(R+1)\ln\left(\frac{R+1}{R}\right)\right]^p=\left[\ln\left(1+\frac{1}{R}\right)^R + \ln\left(1+\frac{1}{R}\right)\right]^p.
\]
Therefore,
\begin{equation}\label{lim5fund.2}
\lim_{R\to +\infty}q(0)=1	
\end{equation}
and
\begin{equation}\label{lim6fund.2}
\lim_{R\to +\infty}q(1)=1,
\end{equation}
where we used the Euler Number 
\[
e=\lim_{R\to +\infty}\left(1+\frac{1}{R}\right)^R.
\]	

By \eqref{lim5fund.2} and \eqref{lim6fund.2}, we obtain that
\[
\lim_{R\to +\infty}q(t)=1 \mbox{ uniformly in } t \in [0,1].
\]
As 
\[
\bar{q}_-(R,R+1) = \int_0^1\min\{1,q(t)\}dt=\int_0^1{\frac{1+q(t) - |1-q(t)|}{2}}dt,
\]
we conclude that
\begin{equation}\label{lim7fund.2}
\lim_{R\to +\infty}\bar{q}_-(R,R+1)=1.
\end{equation}
As 
\[
\bar{q}_+(R,R+1) = \int_0^1\max\{1,q(t)\}dt=\int_0^1{\frac{1+q(t) + |1-q(t)|}{2}}dt,
\]
we conclude that
\begin{equation}\label{lim8fund.2}
\lim_{R\to +\infty}\bar{q}_+(R,R+1)=1.
\end{equation}
It follows from \eqref{estlamb} that
\begin{equation}\label{estlamb3.2}
	\left(\frac{k\pi_p}{\bar{q}_+(R,R+1)}\right)^p\leq \lambda_k(R,R+1) \leq \left(\frac{k\pi_p}{\bar{q}_-(R,R+1)}\right)^p.
\end{equation} 
By \eqref{lim7fund.2}, \eqref{lim8fund.2} and by limits in \eqref{estlamb3.2}, we obtain
\[
\lim_{R\to +\infty}\lambda_k(R,R+1)=(\pi_pk)^p.
\] 

The proof of $(\pi_{p}k)^{p}$ is a solution of \eqref{P7.2.1}, according to Zhang \cite{Zhang} (see also del Pino and Manasevich \cite{dPM}).

\textbf{Proof of (ii):} By \eqref{q}, if $N>p$, we have
\begin{equation}\label{q2}
	q(t)=\left(\frac{p-1}{N-p}\right)^{p}A^{\frac{(p-1)p}{N-p}}(B-t)^{-\frac{p(N-1)}{N-p}}.
\end{equation}
Therefore,
\[
q'(t)=\left(\frac{p-1}{N-p}\right)^{p}A^{\frac{(p-1)p}{N-p}}\frac{p(N-1)}{N-p}(B-t)^{-\frac{p(N-1)}{N-p} -1}>0,
\]
and the function $q$ is increasing at $t$. Hence,
\[
0<q(0)\leq q(t) \leq q(1), ,\,\,\ \forall \,\, t \in [0,1].
\]

Since $p=2$ and $N\geq 3$, by\eqref{q2}, we have
\[
q(t)=(N-2)^{-2}\frac{\left(\frac{(R\overline{R})^{N-2}}{\overline{R}^{N-2} - R^{N-2}}\right)^{\frac{2}{N-2}}}{\left(\frac{(\overline{R})^{N-2}}{\overline{R}^{N-2} - R^{N-2}} - t\right)^{\frac{2(N-1)}{N-2}}}.
\] 
In particular,
\[
q(0)=\frac{1}{(N-2)^2}R^2\left(\frac{\overline{R}^{N-2} - R^{N-2}}{\overline{R}^{N-2}}\right)^2
\]
and
\[
q(1)=\frac{1}{(N-2)^2}\overline{R}^2\left(\frac{\overline{R}^{N-2} - R^{N-2}}{R^{N-2}}\right)^2.
\]
If $R>0$ and $\overline{R}=R+1$, we have
\[
\begin{array}{rcl}
q(0)&=&\frac{1}{(N-2)^2}R^2\left(\frac{(R+1)^{N-2} - R^{N-2}}{(R+1)^{N-2}}\right)^2\\
&=& \frac{1}{(N-2)^2}R^2(R+1-R)^2\left(\frac{(R+1)^{N-3} + (R+1)^{N-4}R+\dots + (R+1)R^{N-4}+R^{N-3}}{(R+1)^{N-2}}\right)^2\\
&=&\displaystyle \frac{1}{(N-2)^2}R^2\left(\frac{\sum_{j=0}^{N-3}{(R+1)^{N-(3+j)}R^j}}{(R+1)^{N-2}}\right)^2\\
\end{array}
\]
and we conclude that
\begin{equation}\label{lim1fund}
\lim_{R\to +\infty}q(0)=1,
\end{equation}
where we used that
\[
R^2\left(\frac{\sum_{j=0}^{N-3}{(R+1)^{N-(3+j)}R^j}}{(R+1)^{N-2}}\right)^2=\left(1-\frac{1}{R+1}\right)^2\left(\sum_{j=0}^{N-3}{\left(1-\frac{1}{R+1}\right)^j}\right)^2 \to (N-2)^2,
\]
as $R \to +\infty$.

Similarly,
\[
\begin{array}{rcl}
q(1)&=&\frac{1}{(N-2)^2}(R+1)^2\left(\frac{(R+1)^{N-2} - R^{N-2}}{R^{N-2}}\right)^2\\
&=& \frac{1}{(N-2)^2}(R+1)^2((R+1)-R)^2\left(\frac{(R+1)^{N-3} + (R+1)^{N-4}R+\dots + (R+1)R^{N-4}+R^{N-3}}{R^{N-2}}\right)^2\\
&=&\displaystyle \frac{1}{(N-2)^2}(R+1)^2\left(\frac{\sum_{j=0}^{N-3}{(R+1)^{N-(3+j)}R^j}}{R^{N-2}}\right)^2\\
\end{array}
\]
and we conclude that
\begin{equation}\label{lim2fund}
\lim_{R\to +\infty}q(1)=1.
\end{equation}

Therefore, by \eqref{lim1fund} and \eqref{lim2fund}, we obtain that
\[
\lim_{R\to +\infty}q(t)=1 \mbox{ uniformly in } t \in [0,1].
\]
As 
\[
\bar{q}_-(R,R+1) = \int_0^1\min\{1,q(t)\}dt=\int_0^1{\frac{1+q(t) - |1-q(t)|}{2}}dt,
\]
we conclude that
\begin{equation}\label{lim3fund}
\lim_{R\to +\infty}\bar{q}_-(R,R+1)=1.
\end{equation}
As 
\[
\bar{q}_+(R,R+1) = \int_0^1\max\{1,q(t)\}dt=\int_0^1{\frac{1+q(t) + |1-q(t)|}{2}}dt,
\]
we conclude that
\begin{equation}\label{lim4fund}
\lim_{R\to +\infty}\bar{q}_+(R,R+1)=1.
\end{equation}
Since by \eqref{pi_p} we have $\pi_2=\pi$, it follows from \eqref{estlamb} that
\begin{equation}\label{estlamb2}
	\left(\frac{k\pi}{\bar{q}_+(R,R+1)}\right)^2\leq \lambda_k(R,R+1) \leq \left(\frac{k\pi}{\bar{q}_-(R,R+1)}\right)^2.
\end{equation} 
By \eqref{lim3fund}, \eqref{lim4fund} and by limits in \eqref{estlamb2}, we obtain
\[
\lim_{R\to +\infty}\lambda_k(R,R+1)=\pi^2k^2.
\]
This proves the item ($ii$). 
\end{proof}


\section{Additional results}

Similar to Theorem \ref{T2}, we can get the following result. Let $r \in \mathbb{N}$. Suppose that $p=r+1$ and $N=2r+1$. Then,
\begin{equation}\label{fundlim3}
	\lim_{R \to +\infty}\bar{q}_+(R,R+1)=\lim_{R \to +\infty}\bar{q}_-(R,R+1)=1.
\end{equation}
In particular, by \eqref{estlamb}
\begin{equation}\label{fundlim4}
	\lim_{R \to +\infty} \lambda_k(R,R+1) = \left(k\pi_p\right)^p=(\pi_{r+1}k)^{r+1}.
\end{equation}

Indeed, since $p=r+1$ and $N=2r+1$, then
\[
\left(\frac{p-1}{N-p}\right)^p=1, \frac{N-p}{p-1}=1, \frac{(p-1)p}{N-p}=r+1, \frac{p(N-1)}{N-p}=2(r+1)
\]
and by \eqref{q2}, we have
\[
q(t)=\frac{\left(\frac{R\overline{R}}{\overline{R} - R}\right)^{r+1}}{\left(\frac{\overline{R}}{\overline{R} - R} - t\right)^{2(r+1)}}.
\] 
If $R>0$ and $\overline{R}=R+1$, we have
\[
q(t)=\frac{(R^2+R)^{r+1}}{(R+1-t)^{2(r+1)}}.
\]
In particular,
\[
q(0)=\left(\frac{R}{R+1}\right)^{r+1}<1
\]
and
\[
q(1)=\left(\frac{R+1}{R}\right)^{r+1}>1.
\]
Therefore,
\begin{equation}\label{lim5fund}
\lim_{R\to +\infty}q(0)=1	
\end{equation}
and
\begin{equation}\label{lim6fund}
\lim_{R\to +\infty}q(1)=1.	
\end{equation}
By \eqref{lim5fund} and \eqref{lim6fund}, we obtain that
\[
\lim_{R\to +\infty}q(t)=1 \mbox{ uniformly in } t \in [0,1].
\]
As 
\[
\bar{q}_-(R,R+1) = \int_0^1\min\{1,q(t)\}dt=\int_0^1{\frac{1+q(t) - |1-q(t)|}{2}}dt,
\]
we conclude that
\begin{equation}\label{lim7fund}
\lim_{R\to +\infty}\bar{q}_-(R,R+1)=1.
\end{equation}
As 
\[
\bar{q}_+(R,R+1) = \int_0^1\max\{1,q(t)\}dt=\int_0^1{\frac{1+q(t) + |1-q(t)|}{2}}dt,
\]
we conclude that
\begin{equation}\label{lim8fund}
\lim_{R\to +\infty}\bar{q}_+(R,R+1)=1.
\end{equation}
It follows from \eqref{estlamb} that
\begin{equation}\label{estlamb3}
	\left(\frac{k\pi_p}{\bar{q}_+(R,R+1)}\right)^p\leq \lambda_k(R,R+1) \leq \left(\frac{k\pi_p}{\bar{q}_-(R,R+1)}\right)^p.
\end{equation} 
By \eqref{lim7fund}, \eqref{lim8fund} and by limits in \eqref{estlamb3} we obtain
\[
\lim_{R\to +\infty}\lambda_k(R,R+1)=(\pi_pk)^p=(\pi_{r+1}k)^{r+1}.
\]

\bibliographystyle{amsplain}

\end{document}